\documentclass[reqno,a4paper,final]{amsart}
\usepackage[left=3cm,right=3cm,top=2cm,bottom=2cm,includeheadfoot]{geometry}
\usepackage{amsmath}
\usepackage{amsfonts}
\usepackage{amssymb}
\usepackage{amsfonts}
\usepackage{MnSymbol}
\usepackage{amsmath}
\usepackage{cancel}
\usepackage{graphicx}
\usepackage[numbers]{natbib}
\usepackage{longtable}
\usepackage{color}

\usepackage[ansinew]{inputenc}
\usepackage{amsthm}
\usepackage[arrow, matrix, curve]{xy}
\usepackage{hyperref}
\hypersetup{%
  pdftitle   = {Computation of generalized Killing spinors on reductive homogeneous spaces},
  pdfkeywords = {supergravity, Killing spinors, isometries,
    homogeneity, supersymmetry, symmetric spaces, M-theory},
  pdfauthor  = {Andree Lischewski},
  pdfcreator = {\LaTeX\ with package \flqq hyperref\frqq}
}
\newtheoremstyle{style}   
  {0.5cm}                 
  {0.5cm}                 
  {}                         
  {}                         
  {\normalfont\bfseries}  
  {\normalfont }{   }
  {}
\theoremstyle{style}

\newtheorem{definition}{Definition}\numberwithin{definition}{section}

\newtheorem{bemerkung}[definition]{Remark}

\newcommand{\ph}{\varphi}
\newcommand{\R}{\mathbb{R}}

\newcommand{\C}{\mathbb{C}}

\newcommand{\g}{\mathfrak{g}}

\definecolor{orange}{rgb}{0.9,0.45,0}
\newcommand{\MUNCH}[1]{\relax}

\allowdisplaybreaks[1]
\begin{document}
\title{Computation of generalized Killing spinors on reductive homogeneous spaces}

\author[Andree Lischewski]{Andree Lischewski}
\address[Andree Lischewski]{Humboldt-Universit\"at zu Berlin, Institut f\"ur Mathematik\\
Rudower Chaussee 25, Room 1.310, D12489 Berlin, Germany}
\email{lischews@mathematik.hu-berlin.de}

\begin{abstract} 
We determine the holonomy of generalized Killing spinor covariant derivatives of the form $D= \nabla + \Omega$ on pseudo-Riemannian reductive homogeneous spaces in a purely algebraic and algorithmic way, where $\Omega : TM \rightarrow \Lambda^*(TM)$ is a left-invariant homomorphism. This is essentially an application of the theory of invariant principal bundle connections defined over homogeneous principal bundles. Moreover, the algorithm allows a computation of the associated Killing superalgebra in certain cases. The procedure is demonstrated by determining the super\texttt{}symmetries of certain homogeneous M2 duals, which arise in M-theory.\\
\newline
\smallskip
\noindent \textbf{\keywordsname} \textit{Killing spinors, homogeneity, holonomy, supergravity}
\end{abstract}
\maketitle

\tableofcontents 

\section{Introduction}  
The holonomy algebra of the Levi Civita connection of a symmetric space $(H/K,g)$ can be computed in purely algebraic terms once the Lie algebraic structure of $\mathfrak{h}$ and $\mathfrak{k}$ is known. In particular, it is easy to determine the space of parallel tensors or spinors on a given symmetric space. More generally, \cite{fu,brm} presents an algorithmic procedure which turns the Killing transport equation on a reductive homogeneous space $H/K$ with $H$-invariant metric $g$ and reductive split $\mathfrak{h} = \mathfrak{k} \oplus \mathfrak{n}$ into something algebraic, making it possible to determine the Lie algebra of Killing vector fields, which in general need not to coincide with $\mathfrak{h}$.\\
Turning to the spinorial analogue, we ask whether also the space of geometric Killing spinors, i.e. solutions $\ph \in \Gamma(S^g)$ of $\nabla^{S^g}_X \ph = \lambda \cdot X \cdot \ph$ for some $\lambda \in \R \cup i \R$ on $(H/K,g)$ can be determined using purely algebraic methods. Killing spinors on pseudo-Riemannian reductive homogeneous spaces have been studied intensively in \cite{bfkg,kath}. The latter reference shows how to construct pseudo-Riemannian reductive homogeneous spaces admitting Killing spinors as so called T-dual spaces of Riemannian homogeneous spaces with Killing spinors, which was one of the first construction methods of Killing spinors in arbitrary signature.\\
More generally than Killing spinors, we consider in this article a pseudo-Riemannian reductive homogeneous space $(M=H/K,g)$ and ask for the existence of spinor fields $\ph \in \Gamma(S^g)$ which are parallel wrt. the connection
\begin{align}
D = \nabla^{S^g} + \Omega \cdot , \label{tre}
\end{align}
where $\nabla^{S^g}$ is the Levi-Civita spin connection, $\Omega: TM \rightarrow Cl(TM,g)$ is $H-$invariant and $\Omega(X)$ acts on spinor fields by Clifford multiplication $\cdot$ for fixed $X \in TM$. Let us elaborate on two examples in which the study of (\ref{tre}) appears naturally.\\
First, generalized Killing spinor equations on Lorentzian manifolds of the form (\ref{tre}) arise in physics in supergravity theories when setting the gravitino variation to zero, i.e. $D-$parallel sections give supersymmetries, cf. \cite{jose2,of12}. In this case $\Omega$ is made up of other bosonic fields of the theory. Often the spinors are additionally subject to further algebraic equations. In this situation one is not only interested in knowing whether there are $D-$parallel spinors but also the dimension of the space of $D-$parallel spinors, i.e. the number of unbroken supersymmetries, is important. For instance, the homogeneity theorem states that backgrounds in certain supergravity theories which preserve a sufficiently large fraction of supersymmetry are automatically homogeneous, cf. \cite{jhom}. Moreover, symmetric backgrounds of 11- and 10-dimensional supergravity have been classified recently in \cite{jose,jn} under the additional assumption that all the other bosonic data are also invariant. However, it is yet unclear which of these solutions to the bosonic field equations are also supersymmetric, i.e. one has to go through the classification list and distinguish those backgrounds which admit solutions to (\ref{tre}). This is work in progress. Group-theoretical methods have earlier been used in \cite{ksh} to construct the Killing spinors of special classes of homogeneous supergravity backgrounds.\\
Second, examples of the generalized Killing spinor equation (\ref{tre}) also appear in a more mathematical context, cf. \cite{am}: Let $(\widehat{N},\widehat{g})$ be a Riemannian spin manifold with oriented hypersurface $(N,g)$. If $\widehat{\ph}$ is a parallel spinor on $\widehat{N}$, then its restriction $\ph$ to $N$ satisfies (\ref{tre}) with $\Omega$ being $\frac{1}{2} \cdot$ the Weingarten tensor of the embedding. Conversely, if $W \in \Gamma(End(TN))$ is any symmetric tensor field and $\ph \in \Gamma(N,S^g)$ is a spinor satisfying (\ref{tre}) with $\Omega = \frac{1}{2} \cdot W$ and all data are real analytic, then there is an ambient space $\widehat{N}= N \times (- \epsilon, \epsilon)$ in which $N$ embeds with Weingarten tensor $W$ and $\ph$ extends to a parallel spinor on $\widehat{N}$.\\
Thus, given a reductive homogeneous spin manifold $(M^n=H/K,g)$ and some $H-$invariant $\Omega \in Sym(TM)$,  one can use (\ref{tre}) to decide whether $(M^n,g)$ can be embedded as oriented hypersurface $M^n \subset \widehat{M}$ with Weingarten tensor $\Omega$ into a space $\widehat{M}$ admitting a parallel spinor.\\
\newline
Motivated by these examples, we ask for an algebraic algorithm which solves (\ref{tre}), i.e. determines the dimension of the space of $D-$parallel spinors on a given pseudo-Riemannian reductive homogeneous space $(H/K,g)$ with reductive split $\mathfrak{h}= \mathfrak{k} \oplus \mathfrak{n}$. Under further generic assumptions, namely that the spin structure is homogeneous, as to be made precise in section \ref{se1}, we find in section \ref{se3} that such an algorithm does indeed exist. The only data which enter are the Lie algebraic structure of $\mathfrak{k} \oplus \mathfrak{n}$ and the $Ad_K$-invariant inner product on $\mathfrak{n}$ which corresponds to $g$. The algorithm is essentially an application of Wang's theorem and the theory of invariant connections on principal bundles as studied in detail in \cite{kn1,kn2}. Thus, the theory which is underlying the algorithm is not new but here it is presented in a way such that it is directly accessible for concrete computations. It has earlier been used in \cite{hammerl1} to compute the conformal holonomy group of the product of two spheres. In section \ref{se5} we show based on our previous results how the Killing superalgebra of a supergravity background defined over a pseudo-Riemannian reductive homogeneous space can be computed purely algebraically. \\
Section \ref{se6} applies these results to an interesting class of 11-dimensional $M-$theory backgrounds studied in \cite{fu}: Motivated by the search for new gravity duals to M2 branes with $N>4$ supersymmetry, equivalently characterized as M-theory backgrounds with Killing superalgebra $\mathfrak{osp}(N|4)$ for $N>4$, one classifies homogeneous M-theory backgrounds with symmetry Lie algebra $\mathfrak{so}(n) \oplus \mathfrak{so}(3,2)$. One finds a number of new backgrounds for $n=5$ of the form $S^4 \times X^7$, where the Lorentzian factor $X^7$ is reductive homogeneous under the action of $SO(3,2)$. However, it remains unclear how much supersymmetries these backgrounds preserve. We study two examples of such backgrounds in section \ref{se6} and determine the space of $D-$parallel spinors using the algorithm developed.\\
The final section \ref{coa} extends the algorithm to conformal geometry and oulines how to solve the twistor equation on a reductive homogeneous space.\\
\newline
\textbf{Acknowledgment} The author gladly acknowledges support from the DFG (SFB 647 - Space Time Matter at Humboldt University Berlin) and the DAAD (Deutscher Akademischer Austauschdienst / German Academic Exchange Service). Furthermore, the author would like to thank José Figueroa-O'Farrill and Noel Hustler for many helpful discussions.
\section{Facts about reductive homogeneous spaces} \label{se1}
Our notation for (reductive) homogeneous spaces follows \cite{aa,hammerl1}: Let $M=H/K$ be a connected homogeneous space for some Lie group $H$ and closed subgroup $K$. We shall in addition assume that $H/K$ is \textit{reductive}, i.e. there exists a -from now on fixed- subspace $\mathfrak{n}$ of $\mathfrak{h}$ such that
\begin{align*}
\mathfrak{h} = \mathfrak{k} \oplus \mathfrak{n} \text{ and } \left[\mathfrak{k},\mathfrak{n} \right] \subset \mathfrak{n}. 
\end{align*}
This allows a natural identification $T_{eK}M \cong \mathfrak{n}$, where $e \in H$ is the neutral element. $H-$invariant tensor fields on $H/K$ correspond to $Ad_{K}:K \rightarrow GL(\mathfrak{n})$- invariant tensors of the same type on $\mathfrak{n}$, where the correspondence is given by evaluating the tensor field at the origin $eK \in H/K$.
Let $g$ be a $H-$\textit{invariant} signature $(p,q)$-metric on $H/K$, i.e. for each $h \in H$ left multiplication $l_h$ with $h$ is an isometry. $g$ corresponds to an $Ad_{K}$-invariant scalar product $\langle \cdot, \cdot \rangle_{\mathfrak{n}}$ of the same signature on $\mathfrak{n}$, i.e. $Ad_{K}$ takes values in $O(\mathfrak{n}, \langle \cdot, \cdot \rangle_{\mathfrak{n}})$.\\
\newline
We briefly describe spin structures on oriented pseudo-Riemannian reductive homogeneous spaces. Consider the $SO(p,q)$-bundle $\mathcal{P}^g \rightarrow M$ of oriented orthonormal frames of $(M,g)$. We have
\begin{align*}
\mathcal{P}^g \cong H \times_{Ad_{K}} SO(\mathfrak{n}),
\end{align*}
and
\begin{align*}
\mathcal{P}^g \times_{SO(p,q)}\R^n \cong TM \cong H\times_{Ad_{K}} \mathfrak{n}, 
\end{align*}
where the latter isomorphism is given by $dl_h d\pi_e X \mapsto [h,X]$ for $X \in T_{eK} \left(H/K \right)$.\\
Any lift of the isotropy representation $Ad_{K}$ to the spin group $Spin(\mathfrak{n})$, i.e. any map $\widetilde{Ad}_{K}:K \rightarrow Spin(\mathfrak{n})$ such that the diagram
\begin{align*}
\begin{xy}
  \xymatrix{
      Spin(\mathfrak{n})    \ar[rd]^{\lambda}    &     \\
      K \ar@{->}[r]^{Ad_{K}} \ar[u]^{\widetilde{Ad}_{K}} &   SO(\mathfrak{n})   \\
  }
\end{xy}
\end{align*}
commutes, allows us to fix a \textit{homogeneous spin structure} $(\mathcal{Q}^g=H \times_{\widetilde{Ad}_{K}}Spin(\mathfrak{n}),f^g)$ of $(M,g)$ (cf. \cite{bfkg}), where $f^g:\mathcal{Q}^g \rightarrow \mathcal{P}^g$ is simply the double covering $\lambda:Spin(p,q) \rightarrow SO(p,q)$ in the second factor. From now on we shall always assume that $(M,g)$ admits a homogeneous spin structure and think of this structure as being fixed. For algebraic properties of Clifford algebras $Cl$, their Clifford groups $Cl^*$ and spinor modules $\Delta$, we refer to \cite{ba81,lm,har}. The real or complex spinor bundle is the associated bundle $S^g := \mathcal{Q}^g \times_{Spin(p,q)}\Delta_{p,q}$.

\section{The algorithm} \label{se3}
Let $(M,g)$ be a connected and oriented $n-$dimensional pseudo-Riemannian reductive homogeneous spin manifold of signature $(p,q)$ with fixed decomposition $M=H/K$, $\mathfrak{h} = \mathfrak{k} \oplus \mathfrak{n}$ and corresponding $Ad_{K}$-invariant inner product $\langle \cdot, \cdot \rangle_{\mathfrak{n}}$ on $\mathfrak{n}$. Let $\nabla=\nabla^{S^g} :\Gamma(S^g) \rightarrow \Gamma(T^*M \otimes S^g)$ denote the spinor covariant derivative induced by the (lift of the) Levi Civita connection.
We want to compute the number of linearly independent spinor fields which are parallel wrt. the \textit{modified covariant derivative}
\begin{align*}
D &= \nabla + \Omega: \Gamma(S^g) \rightarrow \Gamma(T^*M \otimes S^g),\\
D_X \ph &= \nabla_X \ph + \Omega(X) \cdot \ph \text{ for }X \in \mathfrak{X}(M), \ph \in \Gamma(S^g),
\end{align*}
where $\Omega: TM \rightarrow Cl(TM,g) \cong \Lambda^*(TM)$ is a vector bundle homomorphism which is \textit{left-invariant}, i.e. $l_h^* \Omega = \Omega$ for $h \in H$, or in more detail
\begin{align}
l_{h^{-1}}^*\left(\Omega \left(dl_h(X)\right) \right) \stackrel{!}{=} \Omega(X) \text{ for all }X \in TM, h \in H. \label{dufff}
\end{align}

To this end, we introduce the homogeneous $Cl^*(\mathfrak{n},\langle \cdot, \cdot \rangle_{\mathfrak{n}})$-bundle $\overline{\mathcal{Q}}:= H \times_K Cl^*(\mathfrak{n},\langle \cdot, \cdot \rangle_{\mathfrak{n}})$, where $K \rightarrow Cl^*(\mathfrak{n},\langle \cdot, \cdot \rangle_{\mathfrak{n}})$ acts by trivial extension of $\widetilde{Ad}^{H/K}:K \rightarrow Spin(\mathfrak{n}) \subset Cl^*(\mathfrak{n},\langle \cdot, \cdot \rangle_{\mathfrak{n}})$. Obviously, there is a natural $H-$left action $L$ on $\overline{\mathcal{Q}}$ and we call a connection $A \in \Omega^1(\overline{\mathcal{Q}},Cl(\mathfrak{n},\langle \cdot, \cdot \rangle_{\mathfrak{n}}))$ $H-$invariant iff $L_h^* A = A$ for $h \in H$. In the usual way, the covariant derivative $D$ on $S^g \cong \overline{\mathcal{Q}} \times_{Cl^*(\mathfrak{n},\langle \cdot, \cdot \rangle_{\mathfrak{n}})} \Delta_{p,q}$ is induced by a connection $A$ on $\overline{\mathcal{Q}}$.\\
One easily verifies that $D$ is even induced by a $H-$\textit{invariant} connection. Checking this is essentially just a reformulation of the isometry-invariance of the Levi Civita connection and the assumption (\ref{dufff}). However, this observation enables us to use further results from \cite{kn1,kn2} which allow the algebraic computation of the holonomy algebra and curvature of invariant connections defined over homogeneous principal bundles.\\
\newline
Carrying these steps out, results with the mentioned theoretical background in a purely algebraic algorithm: To this end, let
\begin{align*}
T_1,...,T_n \text{ be an oriented orthonormal basis of }\mathfrak{n},\\
L_1,...,L_m \text{ be a basis of }\mathfrak{k}.
\end{align*}
There are constants $d^{ki}_v$ such that $\left[ L_k,T_i\right] = \sum_v d^{ki}_v T_v$.\\
\newline
\textit{Step 1}:\\
The connection $D$ is equivalently encoded (in the sense of \cite{kn1,kn2}) in a linear map 
\begin{align*}
\alpha=\alpha_{g} + \alpha_{\Omega} : \mathfrak{h} \rightarrow Cl(p,q),
\end{align*}
taking values in\footnote{At this point, one has to make a choice: Either, one fixes an orthonormal basis of $\mathfrak{n}$ and works with the Clifford algebra $Cl(p,q)$ of $\R^{p,q}$ only, what we will do, or one works more abstractly in the Clifford algebra $Cl(\mathfrak{n},\langle \cdot, \cdot \rangle_{\mathfrak{n}})$.} the Clifford algebra $Cl(p,q)=Cl(\R^n,\langle \cdot, \cdot \rangle_{p,q})$ of $\R^{p,q}$ with standard pseudo-orthonormal basis $(e_1,...,e_n)$, and it splits into parts 
$\alpha_{g}$ and $\alpha_{\Omega}$, describing $\nabla$ and $\Omega$, respectively. Here, 
$\alpha_{g}: \mathfrak{h} \rightarrow \mathfrak{spin}(p,q) \subset Cl(p,q)$ is given by
\begin{equation} \label{ntensor}
\begin{aligned}
\alpha_g(L_k) &= \frac{1}{2} \cdot \sum_{i<j} d^{ki}_j \cdot e_i \cdot e_j \in \mathfrak{spin}(p,q), \\
\alpha_g(T_i) &= \frac{1}{2} \cdot \sum_{i<j} N^i_{ab} e_a \cdot e_b \in \mathfrak{spin}(p,q),
\end{aligned}
\end{equation}
where (cf. \cite{hammerl1,hammerl2}) $N^i_{ab}=N^i_{[ab]}=\frac{1}{2} \cdot \left(\langle [T_i,T_a]_{|\mathfrak{n}},T_b \rangle_{\mathfrak{n}}- \langle [T_i,T_b]_{|\mathfrak{n}},T_a \rangle_{\mathfrak{n}} - \langle [T_a,T_b]_{|\mathfrak{n}},T_i \rangle_{\mathfrak{n}}\right)$ and for any $A \in \mathfrak{h}$ we let $A_{|\mathfrak{n}}$ denote its projection to $\mathfrak{n}$.

\begin{bemerkung}
More invariantly, the map $\alpha_{g}: \mathfrak{h} \rightarrow \mathfrak{spin}(\mathfrak{n},\langle \cdot, \cdot \rangle_{\mathfrak{n}})$ is given as $\lambda_*^{-1} \circ \widetilde{\alpha}_g$, for $\widetilde{\alpha}_g: \mathfrak{h} \rightarrow \mathfrak{so}(\mathfrak{n},\langle \cdot, \cdot \rangle_{\mathfrak{n}})$ uniquely determined by $\langle \alpha(X)Y,Z \rangle = \frac{1}{2} \left(\langle [X,Y]_{|\mathfrak{n}},Z \rangle_{\mathfrak{n}}- \langle [X,Z]_{|\mathfrak{n}},Y \rangle_{\mathfrak{n}} - \langle [Y,Z]_{|\mathfrak{n}},X \rangle_{\mathfrak{n}} \right)$, where $X \in \mathfrak{h}$ and $Y,Z \in \mathfrak{n}$. 
\end{bemerkung}

The $\Omega-$part $\alpha_{\Omega} : \mathfrak{n} \rightarrow Cl(p,q)$, which lives only on $\mathfrak{n}$, is the evaluation of $\Omega$ at $eK \in M$. More precisely, for fixed $i \in \{1,...,n\}$, we have that 
\[ \Omega(T_i) = \sum_{I} \Omega^i_{I} T_I \in Cl(\mathfrak{n},\langle \cdot, \cdot \rangle_{\mathfrak{n}}). \]
The $\Omega^i_I$ are constants, the sum runs over all multi-indicees $(i_1 < i_2 <...i_k)$ for $k \leq n$ and $T_I:=T_{i_1} \cdot...\cdot T_{i_k}$. In this notation,
\begin{align*}
\alpha_{\Omega}(T_i)=\sum_{I} \Omega^i_{I} e_I \in Cl(p,q).
\end{align*}

\textit{Step 2}:\\
One introduces the curvature map $\kappa:\Lambda^2 \mathfrak{n} \rightarrow Cl(p,q)$, which measures the failure of $\alpha$ being a Lie algebra homomorphism. Concretely, one computes for $i<j$
\begin{align*}
 \kappa(T_i,T_j):= [\alpha(T_i),\alpha(T_j)]_{Cl(p,q)} - \alpha([T_i,T_j]_{\mathfrak{h}}),
\end{align*}
and determines the space $\widehat{Im}(\kappa):= \text{span}\{ \kappa(T_i,T_j) \mid i<j\} \subset Cl(p,q)$ and its dimension.\\
\newline
\textit{Step 3}:\\
The holonomy algebra $\mathfrak{hol}(D) \subset \mathfrak{gl}(\Delta_{p,q}) \cong Cl(p,q)$ (or $Cl(p,q) \cong \mathfrak{gl}(\Delta_{p,q}) \oplus \mathfrak{gl}(\Delta_{p,q})$), where $\Delta_{p,q}$ denotes the appropriate spinor module in signature $(p,q)$, can be determined as follows\footnote{We use the notation from \cite{hammerl1,hammerl2} where this construction is reviewed. Moreover, this reference presents some examples and shows how the procedure can be applied to certain Cartan geometries which allows the determination of the conformal holonomy algebra of conformal structures over homogeneous spaces.}: 
$\mathfrak{hol}(D)$ is the $\mathfrak{h}$-module generated by $\widehat{Im}(\kappa)$, i.e.
\begin{align} \label{sumy}
\mathfrak{hol}(D) = \widehat{Im}(\kappa) + \left[\alpha(\mathfrak{h}),\widehat{Im}(\kappa)\right]_{Cl(p,q)}+\left[\alpha(\mathfrak{h}),\left[\alpha(\mathfrak{h}),\widehat{Im}(\kappa)\right]\right]_{Cl(p,q)}+...,
\end{align}
that is one starts with $\widehat{Im}(\kappa) \subset \mathfrak{hol}(D)$ which has been computed before. One adds all elements of type $[\alpha(T_i \text{ or }L_j),\kappa(T_a,T_b)]$. If dim $\widehat{Im}(\kappa) + \left[\alpha(\mathfrak{h}),\widehat{Im}(\kappa)\right]_{Cl(p,q)}= $dim $\widehat{Im}(\kappa)$, we are already done. Otherwise, one adds elements of $\left[\alpha(\mathfrak{h}),\left[\alpha(\mathfrak{h}),\widehat{Im}(\kappa)\right]\right]_{Cl(p,q)}$ until the dimension of the sum in (\ref{sumy}) becomes stable, which will happen 
after at most $(\text{dim }\Delta_{p,q})^2$ steps.\\
\newline
\textit{Step 4}:\\
Once $\mathfrak{hol}(D)$ is known, one computes its natural action on spinors (obtained by restriction of an irreducible representation of $Cl(p,q)$ on $\Delta_{p,q}$ to $\mathfrak{hol}(D)$. In particular, $Ann(\mathfrak{hol}(D)) := \{ v \in \Delta_{p,q} \mid \mathfrak{hol}(D) \cdot v = 0 \}$ can be computed, which by the holonomy principle is isomorphic to the space of parallel spinors wrt. $D$ (for $M$ being simply-connected).

\begin{bemerkung}
If actually $(M=H/K,g)$ is a symmetric space, the algorithm simplifies: First, the map $\alpha_g$ is simply the trivial extension of $\widetilde{ad}_{K}: \mathfrak{k} \rightarrow \mathfrak{spin}(\mathfrak{n})$, i.e. the tensor $N$ in (\ref{ntensor}) vanishes. Moreover, (\ref{sumy}) simplifies to
\begin{align}
\mathfrak{hol}(D) = \widehat{Im}(\kappa) + \left[\alpha(\mathfrak{n}),\widehat{Im}(\kappa)\right]+\left[\alpha(\mathfrak{n}),\left[\alpha(\mathfrak{n}),\widehat{Im}(\kappa)\right]\right]+... \label{holfo}
\end{align}
Note that (\ref{holfo}) generalizes a well-known formula for the holonomy of the Levi Civita connection on symmetric spaces, i.e.  where $\Omega=0$ and thus also $\alpha(\mathfrak{n}) = 0$.
\end{bemerkung}

\section{The associated Killing superalgebra} \label{se5}

Let us now in addition assume that $(M,g)$ is space-and time oriented, which allows a global symmetric squaring of spinor fields to vector fields, $(\ph_1,\ph_2) \rightarrow V_{\ph_1,\ph_2}$, given by $g(V_{\ph_1,\ph_2},X) = \langle \ph_1, X \cdot \ph_2 \rangle_{S^g}$ for a $Spin^+(p,q)$ invariant inner product on the spinor module\footnote{The definition of $V_{\ph_1,\ph_2}$ might also involve taking the real or imaginary part, depending on $(p,q)$} . Assume moreover that for $\ph_i$ being parallel wrt. $D$, the associated vector is Killing (as true for geometric Killing spinors in certain signatures, cf. \cite{boh}, or Killing spinors in 11-dimensional supergravity).\\
Furthermore we assume that the simply-connected reductive homogeneous space is of the form $(M,g)=(H/K,g)$ such that\footnote{Every $X \in \mathfrak{h}$ generates a Killing vector field $X^*$ on $M$ by setting $X^*(hK):= \frac{d}{dt}_{|t=0} \text{exp}(tX) hK$. We assume that these are all Killing vector fields.} $\mathfrak{h} \cong Kill(M,g)$, the space of all Killing vector fields. Under these assumptions, the determination of the associated \textit{Killing superalgebra} (cf. \cite{jose2}) $\mathfrak{g}=\mathfrak{g}_0 \oplus \mathfrak{g}_1$ is purely algebraic:\\
By definition, the odd part is the space of $D-$parallel spinors, i.e. $\mathfrak{g}_1 \stackrel{\ph \mapsto \ph(eK)}{\cong} \{ v \in \Delta_{p,q} \mid \mathfrak{hol} (D) \cdot v =0 \}$ and the even part is given by Killing vector fields, $\mathfrak{g}_0 = \mathfrak{h} = \mathfrak{k} \oplus \mathfrak{n}$, where we use the isomorphism to $Kill(M,g)$ given by
\begin{align} 
Kill(M,g) \ni X \mapsto \left(\left(\nabla X \right)_{eK}, X(eK) \right).\label{i}
\end{align}
In particular, we identify $\mathfrak{k}$ with a subspace of $\mathfrak{so}(\mathfrak{n}) \cong \mathfrak{spin}(\mathfrak{n})$. \\
The brackets in $\mathfrak{g}=\mathfrak{g}_0 \oplus \mathfrak{g}_1$ can now be computed as follows: The even-even bracket, which is classically the usual Lie bracket on vector fields, is simply (minus) the bracket in $\mathfrak{h}$. The odd-even bracket is classically given as the spinorial Lie derivative
\begin{align*}
L_{X} \ph := \nabla_X \ph + \frac{1}{2} \underbrace{\left(\nabla(X) \right)}_{\in \mathfrak{so}(TM) \cong \Lambda^2(TM)} \cdot \ph \stackrel{D \ph = 0}{=} - \Omega(X) \cdot \ph + \frac{1}{2} {\left(\nabla(X) \right)} \cdot \ph
\end{align*}
for $X \in Kill(M,g)$ and $D \ph = 0$. Thus, in the algebraic picture, for $(\beta,t) \in \mathfrak{k} \oplus \mathfrak{n}=\mathfrak{g}_0$, this corresponds to
\begin{align*}
\g_0 \otimes \g_1 \ni (\beta,t) \otimes v \mapsto -\Omega_{eK}(t) \cdot v + \frac{1}{2} \beta \cdot v \in \g_1,
\end{align*}
In order to express the odd-odd-bracket, which squares a $D-$parallel spinor $\ph$ to its Dirac current $V_{\ph,\ph} \in Kill(M,g)$, we differentiate $V{_{\ph}}:=V_{\ph,\ph}$ to obtain,
\begin{align} \label{frt}
g(\nabla_X V_{\ph},Y) = \langle (\epsilon_1 \cdot Y \cdot \Omega(X) + \epsilon_2 \cdot \Omega(X)^T \cdot Y)\cdot \ph, \ph \rangle_{S^g}, 
\end{align}
where $\epsilon_i$ are $(p,q)$-dependent signs and $\Omega(X)^T$ denotes the transpose of $\Omega(X)$ considered as endomorphism acting on spinors. Thus the bracket is by polarization under the isomorphism (\ref{i}) uniquely determined by
\begin{align*}
S^2 \mathfrak{g}_1 \ni v \circ v \mapsto (\alpha_{v},t_{v}) \in \g_0,
\end{align*}
where $t_v$ is the vector in $\mathfrak{n}$ algebraically associated to $v (\leftrightarrow \ph(eK)) \in \mathfrak{g}_1$ and $\alpha_v \in \mathfrak{so}(\mathfrak{n})$ denotes the skew symmetric endomorphism $\nabla V_{\ph}$ given by (\ref{frt}) evaluated at $eK$.\\
\newline
That is, under the assumptions made the structure of $\mathfrak{g}$ can be calculated in a purely algebraic way and $\mathfrak{g}$ can then be analyzed further via its Levi decomposition.

\section{Application to a class of homogeneous M2-duals} \label{se6}
In general, let $(M,g,F)$ be a classical $M-$theory background, i.e. $(M,g)$ is a 11-dimensional connected Lorentzian spin manifold with mostly $+$metric and Clifford algebra convention $x \cdot x = + ||x||^2$ for $x \in \R^{1,10}$, $F$ is a closed $4-$form and we demand the triple to satisfy the bosonic field equations
\begin{equation} \label{fwedge}
\begin{aligned}
d \ast F &= \frac{1}{2} F \wedge F, \\
Ric(X,Y) &= \frac{1}{2} \langle X \invneg F, Y \invneg F \rangle - \frac{1}{6} g(X,Y) |F|^2.
\end{aligned}
\end{equation}
Setting the gravitino variation to zero in a purely bosonic background yields the Killing spinor equation
\begin{align}
D_X \ph = \nabla_X \ph + \underbrace{\frac{1}{6}(X \invneg F) \cdot \ph + \frac{1}{12} \left(X^{\flat} \wedge F\right) \cdot \ph}_{=:\Omega(X) \cdot \ph}. \label{sugra}
\end{align}
A background of 11-dimensional supergravity is called supersymmetric iff it admits nontrivial solutions to (\ref{sugra}). Motivated by the search for new homogeneous M2 duals, \cite{fu} obtains new families of solutions to (\ref{fwedge}): 
The geometry in this case is locally isometric to a product $M= H/K \cong SO(5)/SO(4) \times SO(3,2)/SO(3)$. Let $L_{ab}, L_{a5}$ denote the standard generators of $\mathfrak{so}(5)$, where $a,b=1,2,3,4$, and let $J_{ij},J_{iA},J_{45}$ denote the standard generators of $\mathfrak{so}(3,2)$, where $i,j=1,2,3$ and $A=4,5$. Then
\begin{align*}
\mathfrak{k}&=\mathfrak{k}_1 \oplus \mathfrak{k}_2 = \text{span} \{L_{ab} \} \oplus \text{span} \{J_{ij} \},\\
\mathfrak{n}&=\mathfrak{n}_1 \oplus \mathfrak{n}_2 = \text{span} \{L_{a5} \} \oplus \text{span} \{J_{iA},J_{45} \} \cong \R^4 \oplus \R^7.
\end{align*}
Given the ordered basis $X_{\mu}=(J_{45},L_{15},...,L_{45},J_{14},...,J_{34},...,J_{15},...,J_{35})$, let $\theta^{\mu}$ denote the canonical dual basis for $\mathfrak{n}^*$. Then an $H-$invariant inner product is given by
\begin{align}
g = -(\theta^0)^2+\gamma_1 ((\theta^1)^2 + ...+(\theta^4)^2) + \gamma_2 ((\theta^5)^2 + ...+(\theta^7)^2)+ \gamma_3 ((\theta^8)^2 + ...+(\theta^{10})^2),\label{metrik}
\end{align}
for real parameters $\gamma_{1,2,3} > 0$. Let $(a_0,...,a_{10})$ denote the $g-$pseudo orthonormal basis of $\mathfrak{n}$ obtained by rescaling elements of $X_{\mu}$ with appropriate positive constants.\\
Special choices of the $\gamma_i$ and specifying certain closed $H-$invariant 4-forms on $M=S^4 \times X^7$ yield $M-$theory backgrounds, i.e. solutions to (\ref{fwedge}). Let us consider preserved supersymmetry of two of them in more detail. To this end, we observe that $H-$invariance of $F$ implies $H-$invariance of $\Omega$ and thus makes the algorithm developed before applicable:\\
\newline
\textbf{A supersymmetric Freund-Rubin background. } In this geometry, $\gamma_1= \frac{4}{9}, \gamma_2 = \gamma_3 = \frac{2}{3}$. The Lorentzian factor $X^7$ admits an invariant Lorentzian Sasaki-Einstein structure. The triple $(M,g,F:= \frac{9}{2} \cdot vol_{S^4})$ is a solution to (\ref{fwedge})  and $F$ is obviously $H-$invariant. We turn to (complex) spinor fields on $(M,g)$: 
Let $(e_0,...,e_{10})$ denote the standard basis of $\R^{1,10}$. By $e_i$ we also label the Clifford action of the vector $e_i$ on the complex spinor module $\Delta_{1,10}^{\C}$. We work with the identification
\[ \Delta_{1,10}^{\C} \cong \Delta_4^{\C} \otimes \Delta_{1,6}^{\C}, \]
(cf. also \cite{ldr}) under which Clifford action becomes
\begin{align}
e_{i=1,...,4} \rightarrow e_i \otimes Id\text{ and }e_{i=0,5,...,10} \rightarrow vol_{S^4} \otimes e_i.
\end{align}
In particular, $vol_{S^4}$ acts as the identity on $\Delta_{1,6}^{\C}$. We fix the isometry $\eta:\mathfrak{n} \rightarrow \R^{1,10}$ mapping $a_i$ to $e_i$. It follows directly from the various definitions that the map $\alpha$ describing $D$ splits into $\alpha= \alpha_1 + \alpha_2$, where $\alpha_i : \mathfrak{h}_i \rightarrow Cl(\mathfrak{n}_i)$. As moreover $\alpha_1(X_1)$ and $\alpha_2(X_2)$ commute when acting on spinors, also $\kappa$ splits. The $S^4$-factor is equipped with a multiple of the round standard metric and therefore symmetric. Whence, $\alpha_1(X) = \frac{1}{6} \cdot \frac{2}{3} \cdot \eta(X) \cdot \eta_*F$ for $X \in \mathfrak{n}_1$. Moreover, using that the adjoint action of $\mathfrak{k}_1$ on $\mathfrak{n}_1$ is identified with the identity map, one obtains for $i \neq j$
\begin{align*}
\kappa_1(L_{i5},L_{j5}) &= [\alpha_1(L_{i5}),\alpha_1(L_{j5})] - \alpha_1([L_{i5},L_{j5}])\\
&=\frac{2}{6^2} \cdot \left( \left(\frac{2}{3}\right)^2 \cdot e_i \cdot e_j \cdot \left(\frac{9}{2} \cdot e_1 \cdot e_2 \cdot e_3 \cdot e_4 \right)^2 \right) - \frac{1}{2} \cdot e_i \cdot e_j =0,
\end{align*}
that is $\kappa_1 \equiv 0$. In particular, $\mathfrak{hol}(\alpha_1) = \{0\}$, which reflects the fact that $S^4$ admits a full space of geometric Killing spinors.
Consequently, 
\begin{align}
\mathfrak{hol}(D) = \mathfrak{hol}(\alpha_2) \subset Cl_{1,6} \subset Cl_{1,10}, \label{stre}
\end{align}
and $\alpha_2$ precisely encodes (in the sense of \cite{kn1,kn2}) the connection $\nabla_X + \frac{3}{8} X \cdot$ on $X^7$. (\ref{stre}) and the holonomy-principle direcly reveal that there is a basis of $D-$parallel spinors of the form $\ph_1 \otimes \ph_2$, where $\ph_1$ is a (combination of) Killing spinors on $S^4$ and $\ph_2$ is a geometric Killing spinor to the Killing number $-\frac{3}{8}$ on $X^7$. On the other hand, we have collected all the algebraic ingredients to compute with the algorithm that the following elements lie in $\mathfrak{hol}(\alpha_2)$:
\begin{equation} \label{eq1}
\begin{aligned}
\kappa_2(a_5,a_{10}) &=e_5 \cdot e_{10} + e_7 \cdot e_8,\\
\kappa_2(a_6,a_{10}) &=e_6 \cdot e_{10} + e_7 \cdot e_9,\\
\kappa_2(a_5,a_9) &=e_5 \cdot e_9 + e_6 \cdot e_8,\\
\kappa_2(a_0,a_{10}) &=e_7  + e_0 \cdot e_{10},\\
\kappa_2(a_0,a_{9}) &=e_6  + e_0 \cdot e_{9}.\\
\end{aligned}
\end{equation}
It is known from the general theory (cf. \cite{boh,kath}) that every Lorentzian Sasaki-Einstein manifold admits 2 geometric Killing spinors, which in our case also follows from running the algorithm from section \ref{se3}. On the other hand, it is easy to verify from the definitions that $\mathfrak{hol}(\alpha_2) \subset \mathfrak{spin}(1,6) \oplus \R^{1,6} \cong \mathfrak{spin}(2,6)$. As $X$ is Lorentzian Einstein Sasaki, we must in fact have, cf. \cite{baer,boh,kath} that $\mathfrak{hol}(\alpha_2) \subset \mathfrak{su}(1,3)$. Moreover, $\lambda_* (\mathfrak{hol}(\alpha_2)) \subset \mathfrak{so}(2,6)$ acts irreducible on $\R^{2,6}$ as follows easily from an inspection of the elements (\ref{eq1})\footnote{An element $h \in \mathfrak{spin}(2,6)$ acts via $\lambda_*$ on $\R^{2,6}$ as $x \mapsto [h,x]=h \cdot x - x \cdot h \in \R^{2,6}$. From this it follows easily that there is no proper subspace of $\R^{2,6}$ preserved by alle elements (\ref{eq1}).}. However, there is no proper subgroup of $SU(1,3)$ which acts irreducile on $\R^{2,6}$ as shown in \cite{scal}. This already implies that 
\[\mathfrak{hol}(\alpha_2)=\mathfrak{su}(1,3), \]
which can also be derived by using the algorithm only. Thus, 
\begin{align*}
\mathfrak{hol}(D) = \mathfrak{su}(1,3) \subset \mathfrak{spin}(2,6) = \mathfrak{spin}(1,6) \oplus \R^{1,6} \subset Cl_{1,6} \subset Cl_{1,10}.
\end{align*}
Thus, $X^7$ is a generic Lorentzian Einstein Sasaki manifold and all $D-$parallel spinors are given by tensor products of geometric Killing spinors on $S^4$ and the 2 linearly independent geometric Killing spinors on $X^7$ which define its Sasaki structure. Real $D-$parallel spinors are obtained by imposing additional Majorana conditions (cf. \cite{har,lm,br}).\\
\newline
\textbf{A circle of backgrounds. } This family $(M,g,F)$ of $M-$theory backgrounds is specified from (\ref{metrik}) by the choice $\gamma_1 = \gamma_2 = \gamma_3 = \frac{4}{9}$ and the $\alpha \in \R$-depended $H-$invariant 4-form
\[ F = F_1 + F_2 -\frac{1}{3} \theta^{1234} + \frac{1}{\sqrt{3}} \theta^0  \wedge \text{Re} \left(e^{i \alpha} (\theta^5+ i \theta^8) \wedge (\theta^6+ i \theta^9) \wedge (\theta^7+ i \theta^{10}) \right). \]
Note that $\alpha$ (and henceforth also $D$) does not split into connections on the factors for this geometry. It is very easy to deduce from the algorithm that there are no $D-$parallel spinors in this situation: In fact, we have with the same notation as in the previous case for $i \neq j$:
\begin{equation} \label{ho1}
\begin{aligned}
\kappa_1(L_{i5},L_{j5}) &= [\alpha(L_{i5}),\alpha(L_{j5})] - \alpha([L_{i5},L_{j5}]) \\
&= \left(\left( -\frac{1}{18} \cdot \frac{4}{9} \cdot \left(\frac{1}{3} \cdot \left( \frac{3}{2}\right)^4 \right) + \frac{2}{{12}^2} \cdot \frac{4}{9} \cdot (\eta_*(F_2))^2  \right) - \frac{1}{2}\right) \cdot e_i \cdot e_j.\\
& \in \mathfrak{hol} (\alpha)
\end{aligned}
\end{equation}
A $D-$parallel spinor requires that $\kappa_1(L_{i5},L_{j5})$ considered as endomorphism acting on spinors has a kernel. As $e_i \cdot e_j$ acts invertible on spinors, the expression in brackets in (\ref{ho1}) must be singular. However, the endomorphism $(\eta_*(F_2))^2$ has eigenvalues $0$ and $\approx \pm \left(\frac{3}{2}\right)^4 \cdot 60,75$ only. Whence (\ref{ho1}) does not annihilate any nonzero spinor, i.e. there are no $D-$parallel spinor fields on $M$. \\
More intuitively, but less rigorous, one sees this by plugging in the ansatz $\ph = \ph_1 \otimes \ph_2$ for a $D-$parallel spinor. $\ph_1$ then has to satisfy a geometric Killing spinor equation on $S^4$ but with the wrong Killing constant. However, note that for the algorithm we did not have to make a particular ansatz for the spinor.

\section{A conformal analogue} \label{coa}
If we turn our attention to conformal (spin) geometry, i.e. metrics which differ by multiplication with a positive function, or conformal supergravity and are interested in the study of first order conformally covariant differential equations on spinors, the geometric Killing spinor equation is replaced by the conformally covariant twistor equation (cf. \cite{bfkg,lei,leihabil,bl})
\begin{align*}
\nabla_X \ph + \frac{1}{n} X \cdot D^g \ph = 0 \text{ for } X \in TM,
\end{align*}
where $D^g$ is the $Spin$ Dirac operator (cf. \cite{ba81,fr}). It's solutions are called twistor spinors. Consider the covariant derivative $\widetilde{D}$ on $S^{g,2}:=S^g \oplus S^g$ given by
\begin{align*}
\widetilde{D}_X \begin{pmatrix} \ph  \\ \phi \end{pmatrix} = \begin{pmatrix} \nabla_X^{S^g} & -X \cdot \\ \frac{1}{2}K^g(X) \cdot & \nabla^{S^g}_X \end{pmatrix} \begin{pmatrix} \ph  \\ \phi \end{pmatrix}=: \nabla_X \begin{pmatrix} \ph  \\ \phi \end{pmatrix} + \widetilde{\Omega}(X) \left(\begin{pmatrix} \ph  \\ \phi \end{pmatrix} \right), 
\end{align*}
where $K^g := \frac{1}{n-2} \cdot \left( \frac{scal^g}{2(n-1)}  \cdot g - Ric^g  \right)$ is the Schouten tensor and $\widetilde{\Omega}:TM \rightarrow End(S^{g,2})$. If $\ph$ is a twistor spinor, then $\begin{pmatrix} \ph  \\ -\frac{1}{n} \cdot D^g \ph \end{pmatrix}$ is $\widetilde{D}$-parallel (cf. \cite{bfkg}), and conversely, the first slot of a $\widetilde{D}$-parallel spinor is a twistor spinor.\\
\newline
Suppose now that $(M=H/K,g)$ is a reductive homogeneous pseudo-Riemannian space with $\mathfrak{h}= \mathfrak{k}\oplus \mathfrak{n}$. We identify $\mathfrak{n} \cong \R^{p,q}$ by means of some fixed orthonormal basis. There is a natural homogeneous bundle over $M$ admitting an invariant connection whose holonomy coincices with that of $\widetilde{D}$. To this end, we enlarge $\R^{p,q}$ to $\R^{p+1,q+1}$ by introducing new lightlike directions $e_{\pm}$ such that $\langle e_+, e_- \rangle = 1$. Clearly, $\R^{p+1,q+1} = \R e_- \oplus \R^{p,q}\oplus \R e_+$ as $O(p,q)-$modules. We define the annihilation spaces $Ann(e_{\pm}):=\{ v \in \Delta_{p+1,q+1} \mid e_{\pm}\cdot v = 0 \}$. It follows that for every $v \in \Delta_{p+1,q+1}$ there is a unique $w \in \Delta_{p+1,q+1}$ such that $v=e_+ w + e_- w$, leading to a decomposition
\begin{align*}
\Delta_{p+1,q+1} = Ann(e_+) \oplus Ann(e_-). 
\end{align*}
$Ann(e_{\pm})$ is acted on by $Spin(p,q) \hookrightarrow Spin(p+1,q+1)$ and there is an isomorphism $\chi: Ann(e_-) \rightarrow \Delta_{p,q}$ of $Spin(p,q)$-modules leading to the identification
\begin{equation}
\begin{aligned} \label{deco}
\Pi: {\Delta_{p+1,q+1}}_{|Spin(p,q)} & \rightarrow \Delta_{p,q} \oplus \Delta_{p,q}, \\
v=e_+w+e_-w & \mapsto (\chi(e_-e_+w),\chi(e_-w)). 
\end{aligned}
\end{equation}
We identify $\Delta_{p+1,q+1}$ and $\Delta_{p,q} \oplus \Delta_{p,q}$ by means of $\Pi$. Now we introduce the bundle 
\[\widehat{Q}:=H \times_{\widetilde{Ad}_K} Spin(p+1,q+1), \] 
where $\widetilde{Ad}_K : K \rightarrow Spin(p,q) \hookrightarrow Spin(p+1,q+1)$. Obviously, $S^{g,2} \cong \widehat{Q} \times_{Spin(p+1,q+1)} \Delta_{p+1,q+1}$ and under this identification we obtain with (\ref{deco}) that
\begin{align*}
\widetilde{D}_X \psi = \nabla \psi + \underbrace{(X \cdot s_+ + \frac{1}{2}K^g(X) \cdot s_-) \cdot \psi}_{=\widetilde{\Omega}(X) \cdot \psi},
\end{align*}
where $\nabla$ is induced by the Levi Civita connection, $s_{\pm}$ are the obvious global lightlike sections in $H \times_{Ad_K} \R^{p+1,q+1}$ and $\psi \in \Gamma(S^{g,2})$. From this description it becomes immediate that $\widetilde{D}$ is induced by a connection $A \in \Omega^1(\widehat{Q},\mathfrak{spin}(p+1,q+1))$ in the usual way. As the curvature tensor is isometry-invariant, $\widetilde{\Omega}$ is $H-$invariant and it follows as for the Killing spinor equation that $A$ is $H-$invariant. That is, $A$ is again equivalently characterized (in the sense of \cite{kn1,kn2}) by a linear map 
\[ \alpha = \alpha_g + \alpha_{\widetilde{\Omega}} : \mathfrak{h} \rightarrow \mathfrak{spin}(p+1,q+1), \]
where $\alpha_g:\mathfrak{h} \rightarrow \mathfrak{spin}(p,q) \hookrightarrow \mathfrak{spin}(p+1,q+1)$ describes the Levi Civita connection and is given by (\ref{ntensor}) and $\alpha_{\widetilde{\Omega}}$ lives only on $\mathfrak{n}$ and is given by
\[ \alpha_{\widetilde{\Omega}}(t) = \widetilde{\Omega}_{ek}(t) = t \cdot e_+ + \frac{1}{2}K^g(t) \cdot e_- \in \mathfrak{spin}(p+1,q+1).\]
The curvature tensor of a reductive homogneous space, and thus also $K^g$, can be computed purely algebraically (cf. \cite{aa}). Again formula (\ref{sumy}) applies and yields $\mathfrak{hol}(\widetilde{D}) \cong \mathfrak{hol}(\alpha) \subset \mathfrak{spin}(p+1,q+1)$. One then computes it's natural action on $\Delta_{p+1,q+1}$ and twistor spinors for $(H/K,g)$ are in bijective correspondence to the trivial subrepresentations of $\mathfrak{hol}(\alpha)$. This yields a purely algebraic procedure to solve the twistor equation on any reductive homogeneous spin manifold.

\small
\bibliographystyle{amsalpha}
\bibliography{literatur}
\end{document}